\documentclass[11pt]{article}
\usepackage[a4paper,margin=28mm]{geometry}
\usepackage[T1]{fontenc}
\usepackage{lmodern,microtype}
\usepackage{xcolor}
\usepackage{amsmath,amssymb,amsthm,amscd,mathtools}
\usepackage{enumitem,booktabs,longtable,array}
\usepackage[colorlinks=true,linkcolor=blue!55!black,citecolor=blue!55!black,urlcolor=blue!55!black]{hyperref}
\numberwithin{equation}{section}
\newtheorem{theorem}{Theorem}[section]
\newtheorem{proposition}[theorem]{Proposition}
\newtheorem{lemma}[theorem]{Lemma}
\newtheorem{corollary}[theorem]{Corollary}

\newtheorem*{theorem*}{Theorem}
\newtheorem*{conjecture*}{Conjecture $\mathrm{SC}_n$}
\theoremstyle{definition}
\newtheorem{definition}[theorem]{Definition}

\theoremstyle{remark}
\newtheorem{remark}[theorem]{Remark}
\newtheorem*{remark*}{Remark}

\newcommand{\Tor}{\operatorname{Tor}}
\newcommand{\Ext}{\operatorname{Ext}}
\newcommand{\Hom}{\operatorname{Hom}}

\newcommand{\pd}{\operatorname{pd}}
\newcommand{\fd}{\operatorname{fd}}
\newcommand{\FP}{\mathsf{FP}}

\newcommand{\id}{\operatorname{id}}

\title{Higher Coherence of Flat Algebras}
\author{Alexander Pavlov \and Mikhail Dusman}
\date{}
\begin{document}
\maketitle
\begin{abstract}
We establish higher-coherence results for flat algebras over commutative
base rings.  If $K$ is noetherian and a $K$-flat algebra $A$ admits an
$A^{\mathrm e}$-resolution of length $n$ with projective term in degree
$n$ and flat terms in lower degrees, then every product of copies of
$A$ has right flat dimension at most $n-1$ and $A$ is left strongly
$n$-coherent.  A projective $n$-th relative bar syzygy gives the
corresponding two-sided conclusion.  If $K$ is coherent and $A$ is
$K$-flat and relatively quasi-free, then $A$ is left and right strongly
$2$-coherent.  In contrast, we construct a countable commutative ring
$k$ that is strongly $n$-coherent for every $n\ge2$, while the
one-variable polynomial algebra $k[z]$ is faithfully flat and relatively
quasi-free over $k$ but fails strong $q$-coherence for every finite
$q$.  Explicit modules witness failure at each positive degree.
These results establish coherence transfer over coherent bases and show
that the analogous higher-degree transfer fails for every $n\ge2$,
even after any fixed finite shift.
\end{abstract}
\begingroup
\makeatletter
\renewcommand{\@makefnmark}{}
\makeatother
\footnotetext{\textit{2020 Mathematics Subject Classification.}
Primary 16P70; Secondary 13E15, 13D02, 16E40.\par
\textit{Key words and phrases.}
Coherent rings, strongly $n$-coherent rings, relatively quasi-free algebras,
FP-injective modules, polynomial algebras, square-zero extensions.}
\endgroup
\section{Introduction}\label{sec:introduction}

Let $K$ be a commutative ring and $A$ a unital $K$-algebra.  Its
enveloping algebra is $A^{\mathrm e}=A\otimes_K A^{\mathrm{op}}$, and
its bimodule of relative differentials is
\[
\Omega^1_{A/K}
=\ker\bigl(A\otimes_K A\xrightarrow{\,\mu\,} A\bigr),
\qquad \mu(a\otimes b)=ab.
\]
Following the relative version of the quasi-free framework of
Cuntz--Quillen~\cite{CuntzQuillen}, the algebra $A$ is
\emph{relatively quasi-free over $K$} if $\Omega^1_{A/K}$ is
projective as a left $A^{\mathrm e}$-module.  The starting point is the
following theorem of Bondal--Zhdanovskiy~\cite{BondalZhdanovskiy}.

\begin{theorem*}
Let $K$ be a commutative noetherian ring and let $A$ be a $K$-flat,
relatively quasi-free $K$-algebra.  Then $A$ is left and right coherent.
\end{theorem*}

For an integer $n\ge0$ and a ring $R$, a module belongs to
$\FP_n(R)$ if it admits a partial finite free resolution through degree
$n$.  It belongs to $\FP_\infty(R)$ if it admits a finite-rank free
resolution in every degree.  The ring $R$ is \emph{strongly left
$n$-coherent} if $\FP_n(R)=\FP_\infty(R)$ for left modules; strong
right $n$-coherence is defined similarly.  We call $R$ \emph{strongly
$n$-coherent} if both conditions hold.  For commutative rings, strong
$0$-coherence is noetherianity and strong $1$-coherence is ordinary
coherence.  These notions suggest the following proposed extension of
the preceding theorem.

\begin{conjecture*}
For every $n\ge0$, if $K$ is a commutative strongly $n$-coherent ring
and $A$ is a $K$-flat, relatively quasi-free $K$-algebra, then $A$ is
left and right strongly $(n+1)$-coherent.
\end{conjecture*}

Our first main result determines the validity of the proposed transfer
for coherent and strongly $n$-coherent bases.

\begin{theorem}\label{thm:intro-main}
The following statements hold.
\begin{enumerate}[label=(\roman*)]
\item Let $K$ be a commutative coherent ring and let $A$ be a
      $K$-flat, relatively quasi-free $K$-algebra.  Then $A$ is left
      and right strongly $2$-coherent.
\item There exists a countable commutative ring $k$ which is strongly
      $n$-coherent for every $n\ge2$, although the polynomial algebra
      $A=k[z]$ is not strongly $q$-coherent for any finite $q$.
      The algebra $A$ is free and faithfully flat over $k$, finitely
      presented as a $k$-algebra, and satisfies
      \[
      \Omega^1_{A/k}\cong A^{\mathrm e},\qquad
      \operatorname{pd}_{A^{\mathrm e}}A=1.
      \]
      For each $q\ge1$ there is an explicitly constructed
      $A$-module in $\FP_q(A)\setminus\FP_{q+1}(A)$.
\end{enumerate}
\end{theorem}

Part~\textup{(i)} is Theorem~\ref{sc:cor:SC1}; it establishes $(SC_1)$
under the original relative quasi-free hypothesis over a coherent
base.  Part~\textup{(ii)} is Theorem~\ref{neg:thm:main}.  Its single
polynomial example disproves $(SC_n)$ for every $n\ge2$ and, more
strongly, every uniform finite coherence-degree shift in that range.

\begin{remark}
The one-step shifts in the noetherian and coherent cases are sharp.
Over a noetherian base, a relatively quasi-free algebra can fail to be
noetherian on both sides.  Under the hypotheses of
Theorem~\ref{thm:intro-main}\emph{(i)}, $A$ can fail to be coherent
on both sides.
\end{remark}

If $K$ is noetherian ($0$-coherent), then we provide two new proofs of 
Bondal--Zhdanovskiy theorem, using product--Tor methods and FP-injective
modules.

Our second main result gives a higher-coherence criterion in terms of
bimodule resolutions over a noetherian base.

\begin{theorem}\label{thm:intro-bimodule}
Let $K$ be commutative noetherian, let $A$ be a $K$-flat
$K$-algebra, and let $n\ge1$.  Suppose there is an exact sequence
of $A^{\mathrm e}$-modules
\[
0\longrightarrow P_n\longrightarrow P_{n-1}
\longrightarrow\cdots\longrightarrow P_0
\longrightarrow A\longrightarrow0
\]
with $P_n$ projective and $P_0,\ldots,P_{n-1}$ flat.
Then, for every set $I$, the right $A$-module $A^I$ satisfies
$\operatorname{fd}_A(A^I)\le n-1$, and $A$ is left strongly
$n$-coherent.
\end{theorem}

Theorem~\ref{thm:intro-bimodule} follows from the product--Tor bound in
Theorem~\ref{sy:thm:product-tor-bound} and the higher-coherence
criterion in Theorem~\ref{sy:thm:abstract-main}.
Theorem~\ref{sy:thm:bar-main} gives its two-sided specialization:
when the $n$-th relative bar syzygy is $A^{\mathrm e}$-projective,
$A$ is left and right strongly $n$-coherent, and products of copies of
$A$ have flat dimension at most $n-1$ on both sides.

The appendix establishes the infinite-variable monomial Betti-number
calculations used in the counterexample.

\section{Preliminaries on \texorpdfstring{$n$}{n}-coherence}\label{sec:prelim}
\subsection{Finite presentations, stabilization, and syzygy comparison}
\label{subsec:finite-presentations}\label{subsec:finite-syzygy-inequalities}

Throughout, rings are associative and unital, modules are unital, and a
$k$-algebra means a unital algebra whose structure map $k\to A$ has central
image.  Unless a side is specified otherwise, modules in a left-handed
argument are left modules.

\begin{definition}\label{sy:def:FPn}
Let $R$ be a ring and let $n\geq0$.  A left $R$-module $M$ is
\emph{finitely $n$-presented} if there is an exact sequence
\[
F_n\longrightarrow F_{n-1}\longrightarrow\cdots
\longrightarrow F_1\longrightarrow F_0\longrightarrow M\longrightarrow0
\]
in which every $F_i$ is a finitely generated free left $R$-module.  We
write $M\in\FP_n(R)$.
\end{definition}

Thus $\FP_0(R)$ consists of finitely generated modules and $\FP_1(R)$ of
finitely presented modules.  The classes form a descending chain
\[
\FP_0(R)\supseteq\FP_1(R)\supseteq\FP_2(R)\supseteq\cdots.
\]
We write $\FP_\infty(R)$ for modules admitting a projective resolution by
finitely generated free modules in every degree.  Equivalently,
$\FP_\infty(R)=\bigcap_{m\geq0}\FP_m(R)$; see, for example,
\cite{BravoPerez}.

\begin{definition}\label{sy:def:ncoh}
For $n\geq0$, the ring $R$ is \emph{left $n$-coherent} if
\[
\FP_n(R)\subseteq\FP_{n+1}(R).
\]
Right $n$-coherence is defined using right modules.  The ring is
$n$-coherent if it is both left and right $n$-coherent.
\end{definition}

We use \emph{strongly left $n$-coherent} as a synonym for
\emph{left $n$-coherent}, and similarly on the right.

For $n=0$ this is left noetherianity, and for $n=1$ it is ordinary left
coherence.  The following theorem identifies the one-step condition with
finite resolutions of arbitrary length.
\begin{theorem}\label{sc:thm:BP-stabilization}
For every ring $R$ and $n\ge0$, the following are equivalent on a fixed
side:
\begin{enumerate}[label=(\roman*)]
\item $R$ is strongly $n$-coherent;
\item $\FP_n(R)$ is thick;
\item $\FP_n(R)=\FP_\infty(R)$.
\end{enumerate}
\end{theorem}

\begin{proof}
This is \cite[Theorem~2.4]{BravoPerez}.  Recall that ``thick'' means closed
under direct summands and satisfying the two-out-of-three property for short
exact sequences.
\end{proof}

\begin{remark}
For a commutative ring $R$, strong $0$-coherence is equivalent to
noetherianity.  Indeed, if every finitely generated module is finitely
presented, then every cyclic module $R/I$ is finitely presented, so every
ideal $I$ is finitely generated.  The converse is the usual noetherian
argument.  Strong $1$-coherence is ordinary coherence.
\end{remark}

The following comparison applies to short exact sequences over arbitrary
associative rings.

\begin{lemma}\label{lem:finite-syzygy-comparison}
Let $B$ be an arbitrary ring and let
\[
  0\longrightarrow L\longrightarrow P\longrightarrow M
  \longrightarrow0
\]
be an exact sequence of left $B$-modules.  For any $m\ge1$, if
$P\in\FP_m(B)$, then
\[
   M\in\FP_m(B)\quad\text{ if and only if }\quad
   L\in\FP_{m-1}(B).
\]
\end{lemma}
\begin{proof}
This is the short-exact-sequence finiteness comparison of
Bravo--P\'erez~\cite[Theorem~1.5, parts (1) and (3)]{BravoPerez},
which records the inequalities for the maximal degrees of finite
presentation of the three modules over an arbitrary associative ring.
\end{proof}

\subsection{Bennis \texorpdfstring{$n$}{n}-flatness and the higher Chase criterion}\label{subsec:bennis-products}

The higher Chase criterion is naturally expressed using Bennis's
notion of $n$-flatness, which tests the top Tor group against finitely
$n$-presented modules.  The definition differs from the similarly
named $\FP_n$-flatness condition.

\begin{definition}\label{sc:def:nflat}
For $n\ge1$, a right $R$-module $P$ is \emph{$n$-flat} in the sense used by
Bennis if
\[
\Tor_n^R(P,F)=0\qquad(F\in\FP_n(R)).
\]
\end{definition}

Bennis $n$-flatness uses $\Tor_n$ against $\FP_n$, whereas
$\FP_n$-flatness in the Bravo--P\'erez convention uses $\Tor_1$
against the same test class.  We use the phrase ``Bennis $n$-flat''
when distinguishing the two notions.

Finite presentation is characterized by the following product--tensor test.

\begin{lemma}
\label{sy:lem:fp-product}
Let $R$ be a ring and let $K$ be a finitely generated left $R$-module.
Assume that for every set $I$ the canonical map
\[
\theta_{I,K}:R^I\otimes_R K\longrightarrow K^I,
\qquad
(r_i)_{i\in I}\otimes x\longmapsto(r_i x)_{i\in I},
\]
is injective.  Then $K$ is finitely presented.
\end{lemma}

\begin{proof}
This is the product--tensor characterization of finite presentation:
see Bennis~\cite[Lemma~2.4, equivalence (1)--(2) with $n=1$]{Bennis},
which applies to arbitrary associative rings.  Since $K$ is finitely
generated, the canonical map $R^I\otimes_R K\to K^I$ is surjective;
the assumed injectivity therefore makes it an isomorphism, as required
by the cited criterion.
\end{proof}

The higher-degree form is given by Bennis's generalized Chase criterion.

\begin{proposition}
\label{sy:prop:tor-implies-ncoh}
Let $R$ be a ring and let $n\geq1$.  Suppose that for every set $I$ and
every $M\in\FP_n(R)$ one has
\begin{equation}
\Tor_n^R(R^I,M)=0,
\label{sy:eq:top-tor-criterion}
\end{equation}
where $R^I$ is the product of copies of the regular right $R$-module.
Then $R$ is left $n$-coherent.
\end{proposition}

\begin{proof}
Apply Bennis~\cite[Theorem~2.6, implication (2)$\Rightarrow$(1)]{Bennis}
with $\mathcal X$ equal to the class of all left $R$-modules.
Under this specialization, an $n$-$\mathcal X$-flat right module is
precisely a right module $P$ satisfying
$\Tor_n^R(P,M)=0$ for every $M\in\FP_n(R)$.
Thus the assumed vanishing for all $R^I$ is exactly condition (2)
in the cited theorem, and its conclusion is left strong
$n$-coherence.
\end{proof}

\subsection{Terminal syzygies and a single-product test}\label{subsec:one-product-general}

A single product of sufficiently large cardinality detects the same
Tor vanishing as all products.

\begin{lemma}\label{lem:terminal-flat-resolution}
Let $B$ be any ring, let $r\ge1$, and let $M$ be a left $B$-module.  Suppose
there is an exact sequence
\[
0\longrightarrow Z_r\xrightarrow{\iota} F_{r-1}
\longrightarrow F_{r-2}\longrightarrow\cdots
\longrightarrow F_0\longrightarrow M\longrightarrow0
\]
where all the $F_j$ are flat over $B$.  Then for every right $B$-module
$P$ there is a natural isomorphism
\begin{equation}\label{eq:terminal-general}
\Tor_r^B(P,M)\simeq
\ker\bigl(P\otimes_BZ_r\xrightarrow{1\otimes\iota}
             P\otimes_BF_{r-1}\bigr).
\end{equation}
For $r=1$, the intermediate terms are absent.
\end{lemma}
\begin{proof}
This is the usual dimension-shifting calculation for Tor from a flat
resolution; see Weibel~\cite[Lemma~3.2.8 and its proof]{Weibel}.
The $r=1$ step identifies $\Tor_1^B(P,M)$ with the kernel of
$P\otimes_B Z_1\to P\otimes_B F_0$, and successive connecting
isomorphisms for the flat terms yield the displayed formula for $r>1$.
\end{proof}

The next lemma reduces injectivity of tensor maps for all products to a
single product of a prescribed cardinality.

\begin{lemma}
\label{lem:coordinate-detection}
Let $B$ be a unital ring, $W$ a right $B$-module, and
$f:Z\to F$ a homomorphism of left $B$-modules.
Fix an infinite cardinal $\kappa\ge |W|$.  The following conditions
are equivalent:
\begin{enumerate}[label=(\roman*)]
\item for every set $I$, the map
$W^I\otimes_B Z\xrightarrow{1\otimes f}W^I\otimes_B F$
is injective;
\item for an index set of cardinality $\kappa$, the same map is injective.
\end{enumerate}
\end{lemma}
\begin{proof}
Only the implication from (ii) to (i) needs an argument.
Suppose that, for some set $I$, a nonzero element
\[
 w=\sum_{\ell=1}^{s}p_\ell\otimes z_\ell
 \in W^I\otimes_B Z
\]
is sent to zero by $1\otimes f$.  The sum is finite, and each
$p_\ell$ is a function $I\to W$.  Define
\[
\phi:I\to W^s,\qquad
\phi(i)=(p_1(i),\ldots,p_s(i)),\qquad
J=\phi(I).
\]
Then $|J|\le |W|^s\le\kappa$.  For every $j\in J$ choose
$i_j\in I$ with $\phi(i_j)=j$.  Define the right $B$-linear maps
\[
 u:W^J\to W^I,\quad u(q)(i)=q(\phi(i)),
 \qquad
 v:W^I\to W^J,\quad v(p)(j)=p(i_j).
\]
They satisfy $vu=1_{W^J}$.  Put $q_\ell(j)=p_\ell(i_j)$.
Then $u(q_\ell)=p_\ell$, so
\[
 w=(u\otimes1_Z)w',\qquad
 w'=\sum_{\ell=1}^s q_\ell\otimes z_\ell
 \in W^J\otimes_B Z.
\]
Because $u$ has a left inverse, $u\otimes1_Z$ is split
injective; hence $w'\ne0$.  The same observation for
$u\otimes1_F$ shows that $(1\otimes f)(w')=0$.

Choose an injection $J\hookrightarrow\kappa$.  Extension by zero
realizes $W^J$ as a direct summand of $W^\kappa$.  Tensoring
preserves this splitting, so the nonzero kernel element $w'$
induces a nonzero element in the kernel of
$W^\kappa\otimes_B Z\to W^\kappa\otimes_B F$.
This contradicts (ii).  Therefore (i) holds.
\end{proof}

\begin{theorem}\label{thm:one-product-top-Tor}
Let $B$ be any ring and $W$ a right $B$-module.  Fix a cardinal
$\kappa\ge\max\{\aleph_0,|W|\}$.  Let $r\ge1$, and let $M$ admit a
length-$r$ flat resolution segment as in
Lemma~\ref{lem:terminal-flat-resolution}.  The following are equivalent:
\begin{enumerate}[label=(\roman*)]
\item $\Tor_r^B(W^I,M)=0$ for every set $I$;
\item $\Tor_r^B(W^\kappa,M)=0$.
\end{enumerate}
\end{theorem}
\begin{proof}
Apply the terminal-syzygy formula, Lemma~\ref{lem:terminal-flat-resolution},
to $P=W^I$.  Lemma~\ref{lem:coordinate-detection} then shows that the
resulting tensor-kernel vanishes for every $I$ exactly when it
vanishes for $I$ of cardinality $\kappa$.
\end{proof}

\begin{corollary}
\label{cor:one-product-Chase}
Let $R$ be any unital ring, $r\ge1$, and
$\kappa_R=\max\{\aleph_0,|R|\}$.  Then $R$ is left strongly
$r$-coherent if and only if
\[
\Tor_r^R(R^{\kappa_R},M)=0
\qquad\text{for every }M\in\FP_r(R).
\]
The right-handed statement holds for right modules.
\end{corollary}
\begin{proof}
Every $M\in\FP_r(R)$ has a finite free $r$-presentation.  Its
$r$th syzygy embeds into the finite free module in degree $r-1$, so
Lemma~\ref{lem:terminal-flat-resolution} applies with $B=W=R$.
Theorem~\ref{thm:one-product-top-Tor} reduces vanishing for every
$R^I$ to the displayed single product.  Finally, Bennis's higher Chase
criterion \cite[Theorem~2.6]{Bennis} identifies that universal vanishing
with left strong $r$-coherence.  The right-handed proof replaces $R$ by
$R^{\mathrm{op}}$, whose cardinality is the same.
\end{proof}

\subsection{\texorpdfstring{$\FP_n$}{FPn}-injectivity}\label{subsec:FPinjectivity}

There is also an $\Ext$-theoretic characterization of strong
$n$-coherence.  Throughout this subsection, modules are left modules
over an arbitrary associative unital ring $R$.

\begin{definition}
\label{sc:def:FPninj}
For $n\ge0$, or $n=\infty$, put
\[
\mathcal I_n(R)=\bigl\{E\in R\text{-}\mathrm{Mod}:
  \Ext_R^1(F,E)=0\ \text{for every }F\in\FP_n(R)\bigr\}.
\]
An object of $\mathcal I_n(R)$ is called \emph{$\FP_n$-injective}.
\end{definition}

In degree zero, $\mathcal I_0(R)$ consists of the injective modules:
Baer's criterion tests injectivity against the cyclic modules
$R/J$, all of which are finitely generated.  The class
$\mathcal I_1(R)$ consists of the usual \emph{FP-injective}, or
\emph{absolutely pure}, modules.  The class $\mathcal I_\infty(R)$
is often called the class of \emph{absolutely clean} modules.
Since the classes $\FP_n(R)$ decrease with $n$, the classes of
relative injectives increase:
\[
\mathcal I_0(R)\subseteq\mathcal I_1(R)\subseteq
\mathcal I_2(R)\subseteq\cdots\subseteq\mathcal I_\infty(R).
\]
These definitions and conventions agree with
Bravo--P\'erez~\cite[Definition~3.1]{BravoPerez}.  Unlike Bennis's
$n$-flatness in Subsection~\ref{subsec:bennis-products}, the
condition here tests $\Ext^1$ for each value of $n$.

\begin{theorem}
\label{sc:thm:FPn-inj-criteria}
Let $R$ be any associative unital ring and $n\ge1$.  The following
conditions are equivalent for left $R$-modules:
\begin{enumerate}[label=(\roman*)]
\item $R$ is left strongly $n$-coherent;
\item $\mathcal I_n(R)=\mathcal I_{n+1}(R)$;
\item $\mathcal I_n(R)$ is coresolving;
\item for every $F\in\FP_n(R)$ and $E\in\mathcal I_n(R)$,
      one has $\Ext_R^q(F,E)=0$ for every $q\ge2$.
\end{enumerate}
In (iii), \emph{coresolving} means that the class contains all
injective modules and is closed under extensions and cokernels of
monomorphisms between its objects.
\end{theorem}
\begin{proof}
This is Bravo--P\'erez~\cite[Theorem~5.5]{BravoPerez},
specialized to the four conditions displayed above.  The theorem
is valid for arbitrary associative rings, and separately on either
side.
\end{proof}

In particular, taking $n=1$ gives the \emph{degree-two coherence
criterion} used below:
\begin{equation}\label{sc:eq:FPinj-coherence-test}
R\text{ is left coherent}
\quad\text{if and only if}\quad
\Ext_R^2(F,E)=0
\quad\bigl(F\in\FP_1(R),\ E\in\mathcal I_1(R)\bigr).
\end{equation}

\begin{lemma}\label{sc:lem:cosyzygy-staircase}
Let $R$ be any associative ring and $n\ge0$.  Suppose
\[
0\longrightarrow E\longrightarrow I\longrightarrow C
\longrightarrow0
\]
is exact with $I$ injective.  If $E\in\mathcal I_n(R)$, then
$C\in\mathcal I_{n+1}(R)$.  More generally, after $j$ successive
injective cosyzygies, an $\FP_n$-injective module becomes
$\FP_{n+j}$-injective.  If $R$ is strongly $n$-coherent with
$n\ge1$, every one of these cosyzygies is already in
$\mathcal I_n(R)$.
\end{lemma}
\begin{proof}
Take $F\in\FP_{n+1}(R)$ and a short exact sequence
$0\to\Omega F\to P\to F\to0$ with $P$ finite free.
By the finite-syzygy comparison,
Lemma~\ref{lem:finite-syzygy-comparison}, the module $\Omega F$
belongs to $\FP_n(R)$.  The two long exact $\Ext$ sequences give
\[
\Ext_R^1(F,C)\cong\Ext_R^2(F,E)
\cong\Ext_R^1(\Omega F,E)=0.
\]
This proves the first assertion, and iteration proves the second.
For strongly $n$-coherent $R$, the preceding theorem gives
$\mathcal I_n(R)=\mathcal I_{n+j}(R)$ for all $j\ge0$,
which proves the last assertion.
\end{proof}

The staircase gives the inclusion of cosyzygies in
$\mathcal I_{n+1}(R)$ for an arbitrary ring.  Strong $n$-coherence
requires the reverse inclusion
$\mathcal I_{n+1}(R)\subseteq\mathcal I_n(R)$, as characterized by
Theorem~\ref{sc:thm:FPn-inj-criteria}.

\begin{proposition}
\label{sc:prop:FPn-restriction}
Let $K\to A$ be a unital ring homomorphism, with $A$ flat as a
right $K$-module.  For every $r\ge0$, including $r=\infty$,
restriction of scalars sends $\mathcal I_r(A)$ into
$\mathcal I_r(K)$.  If $K$ is left strongly $m$-coherent and
$r\ge m$, then restriction sends $\mathcal I_r(A)$ into
$\mathcal I_m(K)$.
\end{proposition}
\begin{proof}
Let $X\in\FP_r(K)$ and $E\in\mathcal I_r(A)$.
Choose a projective resolution $P_\bullet\to X$ over $K$.
Since $A$ is right $K$-flat, the complex
$A\otimes_KP_\bullet\to A\otimes_KX$ is an exact complex of
projective left $A$-modules.  Degreewise induction--restriction
adjunction therefore induces natural isomorphisms
\[
\Ext_A^q(A\otimes_KX,E)\cong\Ext_K^q(X,E)
\qquad(q\ge0).
\]
A finite $r$-presentation of $X$ remains a finite
$r$-presentation after tensoring with $A$.  Thus
$A\otimes_KX\in\FP_r(A)$ and the isomorphism in degree one
proves $E\in\mathcal I_r(K)$ upon restriction.
If $K$ is strongly $m$-coherent, then
$\FP_m(K)=\FP_\infty(K)=\FP_r(K)$ whenever $r\ge m$,
so $\mathcal I_r(K)=\mathcal I_m(K)$.
\end{proof}

In particular, if $K$ is noetherian and $A$ is right $K$-flat,
every FP-injective $A$-module restricts to an injective $K$-module.

\section{Higher coherence from projective relative bimodule syzygies}\label{sec:syzygies}
\subsection{Relative bimodules and the bar resolution}

Let $k$ be commutative and let $A$ be a $k$-algebra.  We write
\[
A^{\mathrm e}=A\otimes_k A^{\mathrm{op}}
\]
for the enveloping algebra.  A left $A^{\mathrm e}$-module is the same
thing as an $A$-bimodule: $(a\otimes b^{\mathrm{op}})$ acts by left
multiplication by $a$ and right multiplication by $b$.

\subsubsection{The relative bar complex}

For $q\geq0$ put
\[
B_q(A/k)=A^{\otimes_k(q+2)}.
\]
The two outer tensor factors provide the $A$-bimodule structure.  The bar
differential is
\[
d_q(a_0\otimes\cdots\otimes a_{q+1})
 =\sum_{j=0}^{q}(-1)^j
 a_0\otimes\cdots\otimes a_ja_{j+1}\otimes\cdots\otimes a_{q+1},
\]
and the augmentation $B_0(A/k)=A\otimes_kA\to A$ is multiplication.
The augmented sequence
\begin{equation}
\cdots\longrightarrow B_2(A/k)\longrightarrow B_1(A/k)
\longrightarrow B_0(A/k)\longrightarrow A\longrightarrow0
\label{sy:eq:relative-bar}
\end{equation}
is exact.  One may verify exactness directly from the usual contracting
homotopy obtained by inserting the unit $1\in A$.

\begin{definition}\label{sy:def:bar-syzygy}
For $n\geq1$, the $n$-th relative bar syzygy is
\[
\Omega^n_{\mathrm{bar}}(A/k)=
\begin{cases}
\ker(B_0(A/k)\to A),& n=1,\\[1mm]
\ker(B_{n-1}(A/k)\to B_{n-2}(A/k)),& n\geq2.
\end{cases}
\]
Thus $\Omega^1_{\mathrm{bar}}(A/k)$ is the usual bimodule
$\Omega^1_{A/k}$ of relative noncommutative one-forms.
\end{definition}

\begin{lemma}
\label{sy:lem:bar-flat}
Assume that $A$ is flat as a $k$-module.  Then every $B_q(A/k)$ is flat as
a left $A^{\mathrm e}$-module.
\end{lemma}

\begin{proof}
There is an $A^{\mathrm e}$-linear isomorphism
\begin{equation}
B_q(A/k)\cong A^{\mathrm e}\otimes_k A^{\otimes_k q},
\label{sy:eq:bar-induced}
\end{equation}
where $A^{\mathrm e}$ acts on the first factor.  Since $A$ is $k$-flat,
its finite tensor power $A^{\otimes_k q}$ is $k$-flat.

Let $X\to Y$ be a monomorphism of right $A^{\mathrm e}$-modules.  Under
\eqref{sy:eq:bar-induced}, tensoring with $B_q(A/k)$ over $A^{\mathrm e}$
identifies with tensoring the underlying $k$-modules with
$A^{\otimes_k q}$:
\[
X\otimes_{A^{\mathrm e}}B_q(A/k)
\cong X\otimes_k A^{\otimes_k q}.
\]
The latter functor preserves monomorphisms because
$A^{\otimes_k q}$ is $k$-flat.  Hence $B_q(A/k)$ is
$A^{\mathrm e}$-flat.
\end{proof}

Consequently, if $\Omega^n_{\mathrm{bar}}(A/k)$ is projective over
$A^{\mathrm e}$, then the truncated bar sequence
\begin{equation}
0\longrightarrow\Omega^n_{\mathrm{bar}}(A/k)
\longrightarrow B_{n-1}(A/k)\longrightarrow\cdots
\longrightarrow B_0(A/k)\longrightarrow A\longrightarrow0
\label{sy:eq:truncated-bar}
\end{equation}
is an $A^{\mathrm e}$-flat resolution of length $n$ whose leftmost term is
projective.  The argument uses flatness of the intermediate bar terms; their
projectivity is a stronger condition that can require additional
hypotheses, such as suitable $k$-projectivity.

\subsection{Two product lemmas over a noetherian base}

\begin{lemma}
\label{sy:lem:noeth-product}
Let $k$ be a commutative noetherian ring and let $(F_i)_{i\in I}$ be a
family of flat $k$-modules.  Put $F=\prod_{i\in I}F_i$.  Then:
\begin{enumerate}[label=\textup{(\alph*)}]
\item $F$ is flat over $k$;
\item for every $k$-module $M$, the canonical map
\begin{equation}
\lambda_M:F\otimes_k M\longrightarrow
\prod_{i\in I}(F_i\otimes_k M)
\label{sy:eq:lambda-base}
\end{equation}
is injective.
\end{enumerate}
\end{lemma}

\begin{proof}
Part~\textup{(a)} is Chase's product-flatness theorem
\cite[Theorem~2.1]{Chase}; see also
\cite[Algebra, Proposition~10.90.6, Tag~05CZ]{StacksProject}.
For part~\textup{(b)}, combine the noetherian finite-presentation
property \cite[Algebra, Lemma~10.31.4, Tag~00FP]{StacksProject},
the product--tensor isomorphism for finitely presented modules
\cite[Algebra, Proposition~10.89.3, Tag~059K]{StacksProject},
and flatness from part~\textup{(a)}.
\end{proof}

\begin{lemma}
\label{sy:lem:projective-comparison}
Let $k$ be commutative noetherian, let $A$ be $k$-flat, let $M$ be a left
$A$-module, and let $P$ be a projective left $A^{\mathrm e}$-module.  For
every set $I$, the map induced by the coordinate projections $A^I\to A$,
\begin{equation}
\Phi_P:
(A^I\otimes_k M)\otimes_{A^{\mathrm e}}P
\longrightarrow
\prod_{i\in I}
\bigl((A\otimes_kM)\otimes_{A^{\mathrm e}}P\bigr),
\label{sy:eq:PhiP}
\end{equation}
is injective.
\end{lemma}

\begin{proof}
First assume $P$ is free, say $P=(A^{\mathrm e})^{(J)}$.  Then
\eqref{sy:eq:PhiP} becomes
\[
\bigoplus_{j\in J}(A^I\otimes_kM)
\longrightarrow
\prod_{i\in I}\bigoplus_{j\in J}(A\otimes_kM).
\]
It factors as
\begin{equation}
\bigoplus_{j\in J}(A^I\otimes_kM)
\longrightarrow
\bigoplus_{j\in J}\prod_{i\in I}(A\otimes_kM)
\longrightarrow
\prod_{i\in I}\bigoplus_{j\in J}(A\otimes_kM).
\label{sy:eq:free-factorization}
\end{equation}
The first arrow is injective by Lemma~\ref{sy:lem:noeth-product}(b), applied
to the family of flat $k$-modules $A$.  The second arrow is injective for
purely set-theoretic reasons: an element of a direct sum has only finitely
many nonzero $j$-coordinates, and if all of its $i$-coordinate images are
zero, every one of those finitely many entries is zero.  Thus
\eqref{sy:eq:PhiP} is injective for free $P$.

For general projective $P$, choose a split monomorphism
$P\hookrightarrow (A^{\mathrm e})^{(J)}$.  Tensoring with any right
$A^{\mathrm e}$-module preserves a split monomorphism, and products of
split monomorphisms are monomorphisms.  Comparing \eqref{sy:eq:PhiP} for
$P$ with the already proved free case gives injectivity for $P$.
\end{proof}

\subsection{Relative Tor and the two-sided bar complex}

The next lemma compares relative Tor using the two-sided bar complex.

For a right $A$-module $N$ and a left $A$-module $M$, the $k$-module
$N\otimes_kM$ is a right $A^{\mathrm e}$-module via
\begin{equation}
(n\otimes m)\cdot(a\otimes b^{\mathrm{op}})=na\otimes bm.
\label{sy:eq:right-Ae-action}
\end{equation}

\begin{lemma}\label{sy:lem:relative-tor}
Let $k$ be commutative, let $A$ be $k$-flat, let $N$ be a right
$A$-module which is flat over $k$, and let $M$ be a left $A$-module.
Then for every $q\geq0$ there is a natural isomorphism
\begin{equation}
\Tor_q^{A^{\mathrm e}}(N\otimes_kM,A)
\cong
\Tor_q^A(N,M).
\label{sy:eq:relative-tor-iso}
\end{equation}
\end{lemma}

\begin{proof}
Use the relative bar resolution \eqref{sy:eq:relative-bar} of $A$ as a left
$A^{\mathrm e}$-module.  By Lemma~\ref{sy:lem:bar-flat}, it is a flat
resolution.  In degree $q$ there is a natural isomorphism
\begin{equation}
(N\otimes_kM)\otimes_{A^{\mathrm e}}B_q(A/k)
\xrightarrow{\sim}
N\otimes_k A^{\otimes_k q}\otimes_kM,
\label{sy:eq:bar-identification}
\end{equation}
given by
\[
(n\otimes m)\otimes
(a_0\otimes a_1\otimes\cdots\otimes a_q\otimes a_{q+1})
\longmapsto
na_0\otimes a_1\otimes\cdots\otimes a_q\otimes a_{q+1}m.
\]
The balancing relation for $A^{\mathrm e}$ is exactly what is needed for
this formula to be well defined.  Under these isomorphisms the bar
differential becomes the usual two-sided bar differential.

It remains to justify that the resulting complex computes
$\Tor_*^A(N,M)$.  The augmented right $A$-module bar complex
\begin{equation}
\cdots\longrightarrow
N\otimes_kA\otimes_kA\longrightarrow
N\otimes_kA\longrightarrow N\longrightarrow0
\label{sy:eq:right-bar-N}
\end{equation}
is exact; the standard contracting homotopy inserts the unit of $A$.
Its degree-$q$ term is of the form
\[
(N\otimes_kA^{\otimes_k q})\otimes_kA.
\]
Because $N$ and $A$ are $k$-flat, the coefficient
$N\otimes_kA^{\otimes_k q}$ is $k$-flat.  Therefore this degree-$q$ term
is flat as a right $A$-module: tensoring a left $A$-module with it is the
same as first forgetting to $k$ and then tensoring with a $k$-flat module.
Thus \eqref{sy:eq:right-bar-N} is a flat right $A$-resolution of $N$.
Tensoring it over $A$ with $M$ gives exactly the complex on the right of
\eqref{sy:eq:bar-identification}.  Its homology is therefore
$\Tor_*^A(N,M)$, proving \eqref{sy:eq:relative-tor-iso}.
\end{proof}

\subsection{Higher product--Tor vanishing from a bimodule resolution}

For $n\geq1$, consider an exact sequence of $A$-bimodules
\begin{equation}
0\longrightarrow P_n\longrightarrow P_{n-1}\longrightarrow\cdots
\longrightarrow P_0\longrightarrow A\longrightarrow0,
\label{sy:eq:abstract-resolution}
\end{equation}
in which $P_n$ is projective as an $A^{\mathrm e}$-module and
$P_0,\ldots,P_{n-1}$ are flat as $A^{\mathrm e}$-modules.
A projective $A^{\mathrm e}$-resolution of length $n$ satisfies these
hypotheses.  So does the truncated relative bar resolution
\eqref{sy:eq:truncated-bar} whenever its $n$-th relative bar
syzygy is projective, by Lemma~\ref{sy:lem:bar-flat}.

\begin{theorem}\label{sy:thm:product-tor-bound}
Let $k$ be a commutative noetherian ring, let $A$ be a $k$-flat
$k$-algebra, and let $n\geq1$.  Assume that $A$ has an exact $A^{\mathrm e}$-module resolution of
the form \eqref{sy:eq:abstract-resolution}, with $P_n$ projective
and $P_0,\ldots,P_{n-1}$ flat.  Then, for every set $I$,
\begin{equation}
\fd_A(A^I)\leq n-1
\label{sy:eq:fd-product}
\end{equation}
for $A^I$ regarded as a right $A$-module.  Equivalently,
\begin{equation}
\Tor_q^A(A^I,M)=0
\qquad
\text{for every left $A$-module $M$ and every }q\geq n.
\label{sy:eq:tor-product-vanishing}
\end{equation}
\end{theorem}

\begin{proof}
Fix a set $I$ and a left $A$-module $M$, and write
\[
N=A^I.
\]
Since $A$ is $k$-flat and $k$ is noetherian,
Lemma~\ref{sy:lem:noeth-product}(a) shows that $N$ is $k$-flat.  Therefore
Lemma~\ref{sy:lem:relative-tor} applies to $N$.

Tensor the resolution \eqref{sy:eq:abstract-resolution} on the left with the
right $A^{\mathrm e}$-module $N\otimes_kM$.  Since the $P_j$ are
$A^{\mathrm e}$-flat, the homology of the resulting complex computes
\[
\Tor_*^{A^{\mathrm e}}(N\otimes_kM,A),
\]
and by Lemma~\ref{sy:lem:relative-tor} this is naturally
$\Tor_*^A(N,M)$.

Because the resolution has no terms in degrees $>n$, we immediately have
\begin{equation}
\Tor_q^A(N,M)=0\qquad(q>n).
\label{sy:eq:q-greater-n}
\end{equation}
It remains to prove the top-degree vanishing at $q=n$.

For each $i\in I$, let $\pi_i:N=A^I\to A$ be the coordinate projection.
The maps $\pi_i\otimes\id_M$ produce a morphism from the complex
\[
(N\otimes_kM)\otimes_{A^{\mathrm e}}P_\bullet
\]
to the product, over $i\in I$, of the complexes
\[
(A\otimes_kM)\otimes_{A^{\mathrm e}}P_\bullet.
\]
In degree $n$ this morphism is precisely
\begin{equation}
(N\otimes_kM)\otimes_{A^{\mathrm e}}P_n
\longrightarrow
\prod_{i\in I}
\bigl((A\otimes_kM)\otimes_{A^{\mathrm e}}P_n\bigr),
\label{sy:eq:top-comparison}
\end{equation}
which is injective by Lemma~\ref{sy:lem:projective-comparison} because
$P_n$ is projective.

Let $z$ be an $n$-cycle in the source complex.  Since there is no term in
degree $n+1$, showing that $H_n$ is zero is equivalent to showing that
$z=0$.  For each $i$, the image $z_i$ of $z$ is an $n$-cycle in the
$i$-th target complex.  Its homology class belongs to
\[
\Tor_n^{A^{\mathrm e}}(A\otimes_kM,A)
\cong \Tor_n^A(A,M)=0,
\]
where the isomorphism is Lemma~\ref{sy:lem:relative-tor} with $N=A$.
Again there is no incoming differential in degree $n$, so vanishing of the
homology group forces the cycle $z_i$ itself to be zero.  Thus the image
of $z$ under \eqref{sy:eq:top-comparison} is zero.  Injectivity of
\eqref{sy:eq:top-comparison} gives $z=0$.

Hence $\Tor_n^A(N,M)=0$.  Together with
\eqref{sy:eq:q-greater-n}, this proves \eqref{sy:eq:tor-product-vanishing} and
therefore \eqref{sy:eq:fd-product}.
\end{proof}

\subsection{The higher coherence theorem}

\begin{theorem}
\label{sy:thm:abstract-main}
Let $k$ be a commutative noetherian ring, let $A$ be a $k$-flat
$k$-algebra, and let $n\geq1$.  If $A$ has an exact $A^{\mathrm e}$-module resolution of the form
\eqref{sy:eq:abstract-resolution} whose term $P_n$ is projective and
whose terms $P_0,\ldots,P_{n-1}$ are flat, then $A$ is left
$n$-coherent.
\end{theorem}

\begin{proof}
Theorem~\ref{sy:thm:product-tor-bound} gives
\[
\Tor_n^A(A^I,M)=0
\]
for every set $I$ and every left $A$-module $M$, hence in particular for
every $M\in\FP_n(A)$.  Proposition~\ref{sy:prop:tor-implies-ncoh} then gives
left $n$-coherence.
\end{proof}

The following theorem applies the product--Tor estimate to the truncated
relative bar complex.  Its case $n=1$ yields the noetherian-base
coherence conclusion.

\begin{theorem}
\label{sy:thm:bar-main}
Let $k$ be a commutative noetherian ring, let $A$ be a $k$-flat
$k$-algebra, and let $n\geq1$.  If
\[
\Omega^n_{\mathrm{bar}}(A/k)
\]
is projective as a left $A^{\mathrm e}$-module, then $A$ is left and right
$n$-coherent.  In fact, for every set $I$,
\[
\fd_A(A^I)\leq n-1
\]
and, on the opposite side,
\[
\fd_{A^{\mathrm{op}}}({}^IA)\leq n-1.
\]
\end{theorem}

\begin{proof}
By Lemma~\ref{sy:lem:bar-flat}, the truncated bar sequence
\eqref{sy:eq:truncated-bar} is an exact $A^{\mathrm e}$-module
resolution of length $n$ with a projective term in degree $n$ and
flat terms in the lower degrees.  Theorem~\ref{sy:thm:abstract-main} and
Theorem~\ref{sy:thm:product-tor-bound} therefore give the left-handed
statements.

For the right-handed statements, apply the same argument to
$A^{\mathrm{op}}$.  The ring $A^{\mathrm{op}}$ is $k$-flat.  Moreover the
relative bar complex of $A^{\mathrm{op}}$ is obtained from that of $A$ by
reversing the bimodule actions.  This operation is an equivalence between
the bimodule categories for $A$ and $A^{\mathrm{op}}$ and preserves
projective modules.  Hence the $n$-th relative bar syzygy of
$A^{\mathrm{op}}$ is projective whenever that of $A$ is projective.
Applying the left-handed result to $A^{\mathrm{op}}$ proves right
$n$-coherence of $A$ and the corresponding flat-dimension bound.
\end{proof}

\begin{corollary}
\label{sy:cor:hoch-dim}
Let $k$ be commutative noetherian, let $A$ be $k$-flat, and let
$n\geq1$.  If
\[
\pd_{A^{\mathrm e}}A\leq n,
\]
then $A$ is left and right $n$-coherent.  If
$1\leq d=\pd_{A^{\mathrm e}}A<\infty$, then every product of copies of
$A$ has right and left flat dimension at most $d-1$, and $A$ is
$m$-coherent for every $m\geq d$.
\end{corollary}

\begin{proof}
A projective $A^{\mathrm e}$-resolution of $A$ of length at most $n$ may
be padded by zero modules to obtain a length-$n$ $A^{\mathrm e}$-module
resolution whose term in degree $n$ is projective and whose lower
terms are flat.  Apply Theorems~\ref{sy:thm:product-tor-bound} and
\ref{sy:thm:abstract-main}, and then the same argument to $A^{\mathrm{op}}$.
For $m\geq d$, strong $d$-coherence implies
$\FP_d=\FP_\infty$ and hence $\FP_m=\FP_\infty$.
\end{proof}

\section{The sharp boundary of the proposed higher-coherence transfer}\label{sec:sharp-boundary}

\subsection{The case \texorpdfstring{$SC_0$}{SC0} via products}

\begin{theorem}\label{sc:cor:SC0}
Let $K$ be commutative noetherian and let $A$ be $K$-flat and relatively
quasi-free.  Then $A$ is left and right coherent.
\end{theorem}
\begin{proof}
For every set $I$, the $K$-module $A^I$ is flat.  Indeed, the
noetherian product-flatness lemma proved in
Section~\ref{sec:syzygies} applies to products of flat $K$-modules.
Since $\Omega^1_{A/K}$ is projective over $A^e$, the sequence
\[
0\longrightarrow\Omega^1_{A/K}\longrightarrow A\otimes_K A
\longrightarrow A\longrightarrow0
\]
is a bimodule resolution whose leftmost term is projective.  The
relative product--Tor theorem,
Theorem~\ref{sy:thm:product-tor-bound}, in degree one yields
\[
\Tor_1^A(A^I,M)=0
\qquad\text{for every left $A$-module $M$ and every $I$.}
\]
Thus every product of copies of the regular right $A$-module is flat.
Chase's product criterion implies left coherence.  Applying the same
argument to $A^{\mathrm{op}}$ gives right coherence.
\end{proof}

\subsection{The case \texorpdfstring{$SC_0$}{SC0} via FP-injective modules}

\begin{lemma}\label{sc:lem:coinduced-coresolution}
Let $A$ be relatively quasi-free over a commutative ring $K$.
If a left $A$-module $E$ is injective as a $K$-module, then there is an
exact sequence
\[
0\longrightarrow E\longrightarrow\Hom_K(A,E)
\longrightarrow\Hom_A(\Omega^1_{A/K},E)\longrightarrow0
\]
of left $A$-modules in which the two modules on the right are injective.
Consequently $\operatorname{id}_A(E)\le1$.
\end{lemma}
\begin{proof}
The defining bimodule sequence for $\Omega^1_{A/K}$ splits as a
sequence of left $A$-modules, using $a\mapsto a\otimes1$.
Applying $\Hom_A(-,E)$ therefore gives the displayed short exact
sequence, with the left $A$-actions arising from right multiplication on
the arguments.  Evaluation identifies the middle term with the
coinduced module $\Hom_K(A,E)$; its action is
$(a\varphi)(b)=\varphi(ba)$.  By induction--coinduction adjunction,
\[
\Hom_A(M,\Hom_K(A,E))\cong\Hom_K(M,E).
\]
The right-hand functor is exact in $M$, so $\Hom_K(A,E)$ is injective.
Since $\Omega^1_{A/K}$ is an $A^e$-projective direct summand of a free
bimodule, $\Hom_A(\Omega^1_{A/K},E)$ is a direct summand of a product of
copies of $\Hom_A(A^e,E)\cong\Hom_K(A,E)$.  Products and direct summands
of injective modules are injective, giving the claim.
\end{proof}

\begin{theorem}\label{sc:thm:FPinj-SC0}
Under the assumptions of Theorem~\ref{sc:cor:SC0}, $A$ is left and
right coherent.
\end{theorem}
\begin{proof}
Let $E$ be an FP-injective left $A$-module.  By
Proposition~\ref{sc:prop:FPn-restriction}, applied with
$m=0$ and $r=1$, its restriction to the noetherian ring $K$
is injective.  Lemma~\ref{sc:lem:coinduced-coresolution} gives
$\operatorname{id}_A(E)\le1$.  In particular,
$\Ext_A^2(F,E)=0$ for every finitely presented left $A$-module
$F$.  Criterion~\eqref{sc:eq:FPinj-coherence-test} proves left
coherence.  The same argument for $A^{\mathrm{op}}$ proves
right coherence.
\end{proof}

\begin{remark}\label{sc:rem:sharp-SC0}
The coherence-degree shift by one in
Theorems~\ref{sc:cor:SC0} and~\ref{sc:thm:FPinj-SC0} is sharp.
Their hypotheses allow $A$ to fail strong $0$-coherence
(noetherianity) on both sides.
Let $k$ be a field and take the free associative algebra
$A=k\langle x,y\rangle$, one of the examples in
\cite[Section~4]{BondalZhdanovskiy}.  The words in $x,y$ form a
$k$-basis of $A$, so $A$ is $k$-flat.  A derivation from $A$ into any
$A$-bimodule is determined by arbitrary values on $x$ and $y$.
Consequently $\Omega^1_{A/k}$ is free over $A^e$ on $dx,dy$.

For $n\ge0$, put $L_n=\sum_{i=0}^n A yx^i$.  This left ideal is
spanned by words ending in one of $y,yx,\ldots,yx^n$.  The word
$yx^{n+1}$ has none of these suffixes, so it belongs to $L_{n+1}$
and lies outside $L_n$.  Thus $A$ fails left noetherianity.
Word reversal is an anti-automorphism of $A$ and gives a strictly
ascending chain of right ideals.  Hence the coherent algebra $A$
fails strong $0$-coherence on both sides.
\end{remark}

\subsection{The case \texorpdfstring{$SC_1$}{SC1}: a one-step relative flat-dimension shift}

\begin{lemma}\label{sc:lem:relative-resolution}
Let $A$ be relatively quasi-free over $K$.  Every left $A$-module $M$
fits into a natural exact sequence
\[
0\longrightarrow\Omega^1_{A/K}\otimes_A M
\longrightarrow A\otimes_KM\longrightarrow M\longrightarrow0.
\]
The leftmost term is a direct summand of $A\otimes_KM^{(J)}$ for some set
$J$.
\end{lemma}
\begin{proof}
The bimodule sequence
$0\to\Omega^1_{A/K}\to A\otimes_KA\to A\to0$ splits as a sequence of
right $A$-modules by the section $a\mapsto1\otimes a$.
Tensoring over $A$ with $M$ gives the exact sequence.  An $A^e$-projective
module is a direct summand of $(A^e)^{(J)}$; after tensoring over the
right $A$-action with $M$, this yields a direct summand of
$(A^e\otimes_AM)^{(J)}\cong A\otimes_KM^{(J)}$.
\end{proof}

\begin{theorem}\label{sc:thm:fd-shift}
Assume $A$ is $K$-flat and relatively quasi-free.  If $N$ is a right
$A$-module with $\operatorname{fd}_K N\le d$, then
$\operatorname{fd}_A N\le d+1$.  The corresponding assertion holds on
the opposite side.
\end{theorem}
\begin{proof}
Flat induction identifies
$\Tor_i^A(N,A\otimes_K X)$ with $\Tor_i^K(N,X)$, for every $K$-module
$X$.  Indeed, induce a projective $K$-resolution of $X$ and use the
natural tensor adjunction degree by degree.
Thus $\Tor_i^A(N,A\otimes_KM)$ vanishes for $i>d$.  By the preceding
lemma, the same holds for
$\Tor_i^A(N,\Omega^1_{A/K}\otimes_A M)$: this is a summand of the
corresponding Tor group for $A\otimes_KM^{(J)}$, and the latter also
vanishes for $i>d$.  The long exact Tor sequence of the relative
induced resolution shows $\Tor_i^A(N,M)=0$ for $i>d+1$, for arbitrary
$M$.  This proves the flat-dimension bound.
\end{proof}

\begin{theorem}\label{sc:cor:SC1}
If $K$ is commutative coherent and $A$ is $K$-flat and relatively
quasi-free, then $A$ is left and right strongly $2$-coherent.
\end{theorem}
\begin{proof}
By Chase's product theorem, arbitrary products of flat $K$-modules
are flat.  Every product $A^I$ is therefore $K$-flat.
The one-step flat-dimension shift gives
$\operatorname{fd}_A(A^I)\le1$, so
$\Tor_2^A(A^I,F)=0$ for every left $A$-module $F$.
The strong-$2$-coherence criterion of Bennis
\cite[Theorems~2.3 and~2.6]{Bennis} implies that $A$ is left strongly
$2$-coherent.  Apply the same argument to $A^{\mathrm{op}}$.
\end{proof}

\begin{remark}\label{sc:rem:sharp-SC1}
The coherence-degree shift by one in Theorem~\ref{sc:cor:SC1} is
sharp.  Its hypotheses allow $A$ to fail strong $1$-coherence
(ordinary coherence).
Take Soublin's example, recalled in \cite[Section~4]{BondalZhdanovskiy}.
The countable product
$K=\prod_{n\ge1}\mathbb Q[[u,v]]$ is coherent, while $A=K[t]$ fails
coherence.  The algebra $A$ is free over $K$ on the powers of $t$.
Under $A^e\cong K[t,s]$, the multiplication map $A^e\to A$ sends $s$
to $t$ and has kernel $(t-s)K[t,s]$.  Since $s-t$ is monic in $s$,
multiplication by $t-s$ identifies $A^e$ with this kernel.  Thus
$A$ is relatively quasi-free over $K$ and satisfies the hypotheses
of Theorem~\ref{sc:cor:SC1}.  It is strongly $2$-coherent and fails
strong $1$-coherence.
\end{remark}

\subsection{A local square-zero mechanism and faithful recovery of a bad ring}

All rings in the remaining construction are commutative.  The local square-zero construction is related to the trivial-extension
examples of Bakkari--Mahdou~\cite{BakkariMahdou}.
A trivial extension $B\ltimes E$ of a ring $B$ by a $B$-module $E$ is the additive
group $B\oplus E$ with multiplication
$(b,e)(b',e')=(bb',be'+b'e)$.  The ideal $0\ltimes E$ is square-zero.

\begin{lemma}\label{neg:lem:local}
Let $(B,\mathfrak n)$ be local, put $\kappa=B/\mathfrak n$,
and let $E$ be an infinite-dimensional $\kappa$-vector space, viewed
as a $B$-module through $B\to\kappa$.  Then $k=B\ltimes E$ is local and
not coherent, and every $\FP_2(k)$-module is finite free.  Consequently
$k$ is strongly $n$-coherent for every $n\ge2$.
\end{lemma}
\begin{proof}
The maximal ideal of $k$ is $\mathfrak r=\mathfrak n\ltimes E$.
The ideal $J=0\ltimes E$ annihilates $\mathfrak r$.
For $M\in\FP_2(k)$, take a minimal finite free cover
$k^a\twoheadrightarrow M$; its kernel $N\subseteq\mathfrak r^a$ is
finitely presented.  To see the latter for this chosen cover, compare
it by Schanuel's lemma with the finite free cover in any fixed finite
$2$-presentation of $M$.
Take a minimal finite free cover $k^b\twoheadrightarrow N$ with
finitely generated kernel $H\subseteq\mathfrak r^b$.  Since $JN=0$,
we have $E^b=Jk^b\subseteq H$.  The projection
$\mathfrak r^b=\mathfrak n^b\oplus E^b\to E^b$ is $k$-linear, because
$\mathfrak n E=0$.  As $H$ is finitely generated, its image is
finite-dimensional over $\kappa$, but it contains $E^b$.  Therefore
$b=0$, so $N=0$ and $M\cong k^a$.
Finally, for nonzero $e\in E$ the annihilator of $(0,e)$ is
$\mathfrak r$.  Its image under projection $\mathfrak r\to E$ is
infinite-dimensional over $\kappa$, so $\mathfrak r$ is not finitely
generated.  In a coherent ring annihilators of individual elements are
finitely generated, giving noncoherence.
\end{proof}

\begin{lemma}\label{neg:lem:quotient}
Let $R\to T=R/I$ be a surjective homomorphism of commutative rings,
with $T\in\FP_\infty(R)$.  For every $T$-module $M$ and every $j\ge0$,
$M\in\FP_j(T)$ if and only if $M\in\FP_j(R)$.
Consequently strong $j$-coherence of $R$ implies strong $j$-coherence
of $T$.
\end{lemma}
\begin{proof}
Finite generation is equivalent over $T$ and over $R$.
For $j\ge1$, choose a surjection $T^a\twoheadrightarrow M$ with kernel
$N$.  Because $T^a\in\FP_\infty(R)$, the elementary short-exact-sequence
finiteness rules (proved in Section~\ref{sec:prelim}) show that
$M\in\FP_j(R)$ exactly when $N\in\FP_{j-1}(R)$.
Over $T$ the same assertion holds because $T^a$ is finite free.
Induction on $j$ proves the equivalence.
If $R$ is strongly $j$-coherent and $M\in\FP_j(T)$, then
$M\in\FP_\infty(R)$ and hence $M\in\FP_\infty(T)$, proving the last
assertion.  A quotient by a non-zero-divisor satisfies the hypothesis
because it has a two-term finite free resolution over $R$.
\end{proof}

\begin{lemma}\label{neg:lem:descent}
Let $R\to S$ be a faithfully flat homomorphism of commutative rings.
If $S$ is strongly $j$-coherent, then so is $R$, for any $j\ge0$.
\end{lemma}
\begin{proof}
Let $M\in\FP_j(R)$, and choose a finite $j$-presentation
$F_j\to\cdots\to F_0\to M\to0$.  Its $j$-th syzygy
$Z_j=\ker(F_{j-1}\to F_{j-2})$, with $Z_0=M$, is finitely generated.
By flatness, $S\otimes_RM\in\FP_j(S)$; strong $j$-coherence makes
$S\otimes_RZ_j$ finitely presented.  Finite presentation descends under
faithfully flat maps, so $Z_j$ is finitely presented over $R$.
Therefore the original finite partial resolution extends by one finite
free term; $M\in\FP_{j+1}(R)$.  (For $j=0$, use finite generation
descent directly.)
\end{proof}

\begin{lemma}\label{neg:lem:realization}
For any nonzero local commutative ring $(D,\mathfrak m)$, put
$B=D[t]_{(\mathfrak m,t)}$, $E=(B/(\mathfrak m,t)B)^{(\mathbb N)}$,
$k=B\ltimes E$, and $s=(t,0)\in k$.  Then
\[
C=k[s^{-1}]\cong B[t^{-1}],
\]
and $D\to C$ is faithfully flat.  Every $\FP_2(k)$-module is finite
free.  If $D$ is not strongly $q$-coherent for some $q\ge0$, then
$C$ is not strongly $q$-coherent.
\end{lemma}
\begin{proof}
Since $tE=0$, localization kills $E$ and identifies $C$ with
$B[t^{-1}]$.  The map $D\to C$ is flat because it is a polynomial
extension followed by localizations.  Every prime
$\mathfrak p\subseteq D$ lifts to
$\mathfrak pD[t]\subseteq(\mathfrak m,t)$, and this prime avoids $t$.
Hence it survives both localizations and contracts to $\mathfrak p$.
The induced map $\operatorname{Spec}C\to\operatorname{Spec}D$ is
surjective, which gives faithful flatness.
The other assertions follow from Lemmas~\ref{neg:lem:local} and
\ref{neg:lem:descent}.
\end{proof}

\subsection{The counterexample: one polynomial variable defeats every finite shift}

\begin{theorem}\label{neg:thm:main}
There exist a countable commutative ring $k$ and a countable,
faithfully $k$-flat, relatively quasi-free $k$-algebra $A=k[z]$ such that
\begin{enumerate}[label=(\roman*)]
\item every module in $\FP_2(k)$ is finite free, so $k$ is strongly
$n$-coherent for all $n\ge2$, but $k$ is not coherent;
\item $A$ is a polynomial algebra in one variable, free as a $k$-module
and finitely presented as a $k$-algebra;
\item $\Omega^1_{A/k}\cong A^e$ and $\operatorname{pd}_{A^e}A=1$;
\item $A$ is not strongly $q$-coherent for \emph{any} finite $q$.
\end{enumerate}
More precisely, for each $q\ge1$ there exists an explicitly defined
$A$-module in $\FP_q(A)\setminus\FP_{q+1}(A)$.
\end{theorem}
\begin{proof}
We use the quadratic monomial ring of Bravo--P\'erez~\cite[Example~1.3]{BravoPerez}.
Appendix~\ref{app:sec:Quadratic monomial algebras and their quotients}
computes its syzygies by a monomial Tor recursion.
Fix $\mathbb F=\mathbb Q$ and define
\[
R=\frac{\mathbb F[x_1,x_2,\ldots,y_1,y_2,\ldots]}
{(x_jx_{j+1},\ x_1y_1,\ y_1y_i\ :\ i,j\ge1)},
\qquad \mathfrak m=(x_i,y_i:i\ge1),\quad D=R_{\mathfrak m}.
\]
Theorem~\ref{app:thm:FP} gives
\[
N_q=(x_q)D\in\FP_q(D)\setminus\FP_{q+1}(D)
\qquad(q\ge1).
\]
In particular $D$ fails strong $q$-coherence for all finite $q\ge1$.
Form $B=D[t]_{(\mathfrak m D,t)}$, put
$E=\mathbb F^{(\mathbb N)}$ with $B$ acting via its residue field,
$k=B\ltimes E$, $s=(t,0)$, and $C=k[s^{-1}]$.
Lemma~\ref{neg:lem:local} proves that $k$ is strongly $n$-coherent
for all $n\ge2$.  By Lemma~\ref{neg:lem:realization},
$D\to C$ is faithfully flat and $C=B[t^{-1}]$.

Set $A=k[z]$.  It is free and faithfully flat over $k$.
The enveloping algebra is $A^e\cong k[z,w]$; multiplication has kernel
$(z-w)k[z,w]$.  Since the monic polynomial $z-w$ is a non-zero-divisor
in $k[z,w]$, multiplication gives
$\Omega^1_{A/k}\cong A^e$.  The two-term free bimodule resolution
of $A$ shows $\operatorname{pd}_{A^e}A\le1$, and equality holds as follows.  If $A$ were projective over
$P=A^e=k[z,w]$, the kernel $I=(z-w)P$ of the surjection $P\to A$
would be a direct-summand ideal, hence $I=eP$ for an idempotent
$e\in P$.  But $(1-e)(z-w)=0$ and $z-w$ is regular, so $e=1$.
This would give $I=P$, contradicting the nonzero quotient $A$.
Therefore $A$ is not projective over $A^e$ and its projective
dimension is exactly one.

Let $f=sz-1\in A$.  It is a non-zero-divisor even though $s$ is a
zero divisor: if $f\sum_{i=0}^m a_iz^i=0$, coefficient comparison
first gives $a_0=0$ and then $a_i=sa_{i-1}=0$ successively.
The quotient $A/(f)$ is $k[s^{-1}]\cong C$; it lies in
$\FP_\infty(A)$ because $f$ is regular.
If $A$ were strongly $q$-coherent for some finite $q$,
Lemma~\ref{neg:lem:quotient} would make $C$ strongly
$q$-coherent.  Then the faithfully flat descent
Lemma~\ref{neg:lem:descent} would make $D$ strongly
$q$-coherent, contradicting the preceding explicit modules $N_q$.

For each $q\ge1$, put
$M_q=C\otimes_D(x_q)D$, regarded as an $A$-module through
$A\twoheadrightarrow C$.  Faithful flatness preserves their
finite $q$-presentations and reflects the failure of
$(q+1)$-presentations (apply finite-presentation descent to the
terminal syzygy).  Thus
$M_q\in\FP_q(C)\setminus\FP_{q+1}(C)$, and
Lemma~\ref{neg:lem:quotient} gives
$M_q\in\FP_q(A)\setminus\FP_{q+1}(A)$.
Countability follows from countability of the variables, the field
$\mathbb Q$, polynomial algebras, quotients, finite localizations,
and the countable direct sum $E$.
\end{proof}

\appendix
\section{Quadratic monomial algebras: linearity, Betti numbers, and syzygies}\label{app:sec:Quadratic monomial algebras and their quotients}
Let $\mathbb F$ be a field, let $V$ be any set, and let $\Gamma$ be an undirected graph on $V$ in which loops are allowed.  Write $u\sim v$ if $\{u,v\}$ is an edge, with $u=v$ allowed.  Let
\[
R=\mathbb F[V]/(uv:u\sim v),\qquad R_+=\bigoplus_{j\ge1}R_j,
\]
where every variable has degree one.  Thus $R$ is a connected nonnegatively graded commutative $\mathbb F$-algebra with $R_0=\mathbb F$; it may have infinitely many generators.  In the expressions above, the notation $\mathbb F[V]$ means the polynomial algebra on the set $V$, i.e. the direct union of its finite-variable subalgebras.  For $S\subseteq V$, write
\[
I_S=(v:v\in S)\subseteq R,\qquad Q_S=R/I_S,\qquad N(v)=\{w\in V:w\sim v\}.
\]
The set $N(v)$ may be infinite.  We regard $\mathbb F=R/R_+$ as a right $R$-module via the augmentation.  A graded shift is denoted $M(-1)$, so that $(M(-1))_d=M_{d-1}$.

For a graded left $R$-module $M$, its \emph{graded Betti numbers} and \emph{total Betti numbers} over $R$ are
\[
\beta_{i,j}^{R}(M)
 =\dim_{\mathbb F}\Tor_i^R(\mathbb F,M)_j,
\qquad
\beta_i^{R}(M)
 =\dim_{\mathbb F}\Tor_i^R(\mathbb F,M)
\quad(i\ge0).
\]
These dimensions may be infinite cardinals.

\begin{lemma}\label{app:lem:colon}
For every $S\subseteq V$ and $v\in S$, if $J=S\setminus\{v\}$ and $T=J\cup N(v)$, multiplication by $v$ gives an exact sequence of graded $R$-modules
\[
0\longrightarrow Q_T(-1)\xrightarrow{\ \cdot v\ }Q_J\longrightarrow Q_S\longrightarrow0.
\tag{1}\label{app:eq:colon-sequence}
\]
\end{lemma}
\begin{proof}
The monomials not divisible by an edge monomial $uw$ form a $\mathbb F$-basis of $R$.  In $Q_J$, the surviving monomials are exactly the admissible monomials not containing any variable from $J$.  Multiplication by $v$ takes such a monomial to zero precisely when it contains a factor $w\in N(v)$; when the product is nonzero, distinct monomials have distinct products.  Hence the annihilator of $v$ in $Q_J$ is generated by the images of $N(v)$, or equivalently
\[
(I_J:_R v)=I_{J\cup N(v)}=I_T.
\]
The image of multiplication by $v$ is $vQ_J$, and its cokernel is $Q_J/vQ_J=Q_S$.  The shift makes the multiplication map homogeneous of degree zero.  These observations give \eqref{app:eq:colon-sequence}.
\end{proof}

\begin{lemma}\label{app:lem:finite}
If $V$ is finite, then for every $S\subseteq V$ and every $i\ge0$,
\[
\Tor_i^R(\mathbb F,Q_S)_j=0\qquad\text{whenever }j\ne i.
\tag{2}\label{app:eq:finite-linearity}
\]
\end{lemma}
\begin{proof}
We use induction on the homological degree $i$.  For $i=0$,
$\Tor_0^R(\mathbb F,Q_S)=\mathbb F\otimes_RQ_S=\mathbb F$, concentrated in degree zero.

Fix $i\ge1$, and assume the assertion is true in degree $i-1$ for every variable subset.  Prove it in degree $i$ by induction on $|S|$.  For $S=\varnothing$ one has $Q_S=R$, so $\Tor_i^R(\mathbb F,Q_S)=0$.  For nonempty $S$, choose $v\in S$ and use the notation $J,T$ of Lemma~\ref{app:lem:colon}.  The long exact graded Tor sequence for \eqref{app:eq:colon-sequence} contains
\[
\Tor_i^R(\mathbb F,Q_J)\longrightarrow
\Tor_i^R(\mathbb F,Q_S)\longrightarrow
\Tor_{i-1}^R(\mathbb F,Q_T(-1)).
\]
The first term is concentrated in internal degree $i$ by induction on $|S|$.  The third is concentrated in degree $i$, because
\[
\Tor_{i-1}^R(\mathbb F,Q_T(-1))
\cong\Tor_{i-1}^R(\mathbb F,Q_T)(-1)
\]
and the outer induction puts $\Tor_{i-1}^R(\mathbb F,Q_T)$ in degree $i-1$.  Thus the middle term has zero homogeneous component outside degree $i$.  This completes both inductions.
\end{proof}

\begin{lemma}\label{app:lem:bar}
Let $H=\mathbb F\oplus H_+$ be any augmented associative $\mathbb F$-algebra.  For $i\ge0$ put
\[
P_i(H)=H_+^{\otimes_\mathbb F i}\otimes_\mathbb F H.
\]
These are free right $H$-modules.  With $P_0(H)=H$, define $d_0=\epsilon:H\to\mathbb F$ and, for $i\ge1$,
\begin{align*}
d_i([a_1|\cdots|a_i]\otimes h)
={}&\sum_{j=1}^{i-1}(-1)^{j-1}
[a_1|\cdots|a_ja_{j+1}|\cdots|a_i]\otimes h\\
&+(-1)^{i-1}[a_1|\cdots|a_{i-1}]\otimes a_ih.
\end{align*}
Then $P_\bullet(H)\to\mathbb F$ is a free right $H$-resolution of the augmentation module $\mathbb F$.
\end{lemma}
\begin{proof}
This is the normalized bar resolution of the augmentation module,
written for right modules.  It follows from Weibel
\cite[Section~8.6.12 and Exercise~8.6.4]{Weibel},
applied to the opposite algebra and specialized to the augmentation
$H\to\mathbb F$.
\end{proof}

\begin{lemma}\label{app:lem:colimit}
Let $R=\mathbb F[V]/(uv:u\sim v)$ with $V$ arbitrary and $S\subseteq V$.  For a finite subset $U\subseteq V$, set
\[
R_U=\mathbb F[U]/(uv:u\sim v,\ u,v\in U),
\qquad Q_{S,U}=R_U/(S\cap U).
\]
Then, for each $i\ge0$, there is a natural isomorphism of graded $\mathbb F$-vector spaces
\[
\Tor_i^R(\mathbb F,Q_S)
\cong\varinjlim_{U\subseteq V,\ U\text{ finite}}
\Tor_i^{R_U}(\mathbb F,Q_{S,U}).
\tag{3}\label{app:eq:tor-colimit}
\]
\end{lemma}
\begin{proof}
The admissible monomial bases show that the transition maps $R_U\to R_W$ and $Q_{S,U}\to Q_{S,W}$ for $U\subseteq W$ are injective.  They also show
\[
R=\bigcup_U R_U,
\quad R_+=\bigcup_U(R_U)_+,
\quad Q_S=\bigcup_U Q_{S,U}.
\]
Tensor Lemma~\ref{app:lem:bar}'s resolution with $Q_S$.  The resulting chain complex computing $\Tor_i^R(\mathbb F,Q_S)$ has terms
\[
C_i(R,Q_S)=R_+^{\otimes_\mathbb F i}\otimes_\mathbb F Q_S.
\]
A typical chain is a finite sum of pure tensors.  Every factor in every such tensor is a finite $\mathbb F$-linear combination of monomials and hence involves only finitely many variables.  Therefore the chain already lies in
\[
C_i(R_U,Q_{S,U})=(R_U)_+^{\otimes_\mathbb F i}\otimes_\mathbb F Q_{S,U}
\]
for some finite $U$.  The bar differential uses only multiplication, which commutes with all inclusions.  Thus
\[
C_\bullet(R,Q_S)=\varinjlim_U C_\bullet(R_U,Q_{S,U})
\]
as graded complexes of $\mathbb F$-vector spaces.

Homology commutes with this filtered union: if a cycle becomes a boundary in the union, it is the boundary of a chain involving finitely many variables, so the bounding chain and the cycle both lie at a common finite stage.  Similarly, each cycle in the union occurs at a finite stage.  Hence cycles modulo boundaries commute with the union, in each internal degree.  Using Lemma~\ref{app:lem:bar} at every stage proves \eqref{app:eq:tor-colimit}.
\end{proof}

\begin{theorem}\label{app:thm:linearity}
For an arbitrary (possibly infinite) variable set $V$, every variable subset $S\subseteq V$ and every $i\ge0$,
\[
\Tor_i^R(\mathbb F,Q_S)_j=0\qquad(j\ne i).
\tag{4}\label{app:eq:linearity}
\]
\end{theorem}
\begin{proof}
For each finite $U\subseteq V$, Lemma~\ref{app:lem:finite} applied to $R_U$ and $S\cap U$ gives concentration of $\Tor_i^{R_U}(\mathbb F,Q_{S,U})$ in degree $i$.  Equation \eqref{app:eq:tor-colimit} is an isomorphism of graded vector spaces, so the same concentration holds after the filtered union.
\end{proof}

\begin{theorem}\label{app:thm:rec}
Let $S\subseteq V$ and consider the graded $R$-module $Q_S=R/I_S$.
Its first two total Betti numbers are
\[
\beta_0^R(Q_S)=1,
\qquad
\beta_1^R(Q_S)=|S|.
\tag{5}\label{app:eq:initial-betti}
\]
If $v\in S$, set $J=S\setminus\{v\}$ and $T=J\cup N(v)$.
For every $i\ge1$ there is a short exact sequence of graded
$\mathbb F$-vector spaces
\[
0\longrightarrow\Tor_i^R(\mathbb F,Q_J)
\longrightarrow\Tor_i^R(\mathbb F,Q_S)
\longrightarrow\Tor_{i-1}^R(\mathbb F,Q_T)(-1)
\longrightarrow0.
\tag{6}\label{app:eq:tor-sequence}
\]
Consequently, the total Betti numbers of these three $R$-modules satisfy
\[
\beta_i^R(Q_S)
 =\beta_i^R(Q_J)+\beta_{i-1}^R(Q_T)
\qquad(i\ge1).
\tag{7}\label{app:eq:betti-recurrence}
\]
The sum in \eqref{app:eq:betti-recurrence} is cardinal addition; the groups and dimensions may
be infinite.
\end{theorem}
\begin{proof}
We first calculate the initial Betti numbers of the module $Q_S$.
Since $I_S\subseteq R_+$, the augmentation gives
$\mathbb F\otimes_R Q_S\cong\mathbb F$.  Thus
$\beta_0^R(Q_S)=1$.

Apply $\mathbb F\otimes_R-$ to the defining exact sequence
\[
0\longrightarrow I_S\longrightarrow R\longrightarrow Q_S
\longrightarrow0.
\]
Since $R$ is free over itself and the map
$\mathbb F\otimes_R I_S\to\mathbb F\otimes_R R$ is zero,
the associated long exact Tor sequence gives
\[
\Tor_1^R(\mathbb F,Q_S)
 \cong\mathbb F\otimes_R I_S
 \cong I_S/R_+I_S.
\]
The ideal $I_S$ is generated by the degree-one variables in $S$.
The classes of these variables form an $\mathbb F$-basis of
$I_S/R_+I_S$: they span because every other monomial in $I_S$
has a factor in $R_+$, and they are independent because the defining
monomial relations are quadratic.  Hence
$\Tor_1^R(\mathbb F,Q_S)\cong\mathbb F^{(S)}(-1)$
and $\beta_1^R(Q_S)=|S|$.

Now fix $v\in S$.  By Lemma~\ref{app:lem:colon},
\[
0\longrightarrow Q_T(-1)\longrightarrow Q_J
\longrightarrow Q_S\longrightarrow0
\]
is exact.  Its long exact graded Tor sequence contains
\[
\begin{aligned}
\Tor_i^R(\mathbb F,Q_T)(-1)
 &\longrightarrow\Tor_i^R(\mathbb F,Q_J)
 \longrightarrow\Tor_i^R(\mathbb F,Q_S)\\
 &\longrightarrow\Tor_{i-1}^R(\mathbb F,Q_T)(-1)
 \longrightarrow\Tor_{i-1}^R(\mathbb F,Q_J).
\end{aligned}
\]
Theorem~\ref{app:thm:linearity} places the unshifted group
$\Tor_h^R(\mathbb F,Q_U)$ in internal degree $h$, for every
variable subset $U$.  The two outer maps displayed above therefore
vanish: in each case their source is concentrated one internal degree
higher than their target.  Exactness gives \eqref{app:eq:tor-sequence}.
Taking $\mathbb F$-dimensions in \eqref{app:eq:tor-sequence} yields \eqref{app:eq:betti-recurrence}.  This remains
valid for infinite-dimensional vector spaces, since an exact
sequence of vector spaces splits and dimensions add as cardinals.
\end{proof}

For a set of vertices $S$, let $\mathcal B_d(S)$ denote the set of vertices at distance at most $d$ from $S$ in $\Gamma$.  A loop does not change distances.  By convention $\mathcal B_d(\varnothing)=\varnothing$.  The following lemma uses this graph-theoretic condition to bound the Betti numbers of the module $Q_S$.
\begin{lemma}\label{app:lem:radius}
Let $S\subseteq V$ be finite and $i\ge1$.  If $\mathcal B_{i-1}(S)$ is finite, then $\beta_i^R(Q_S)$ is finite.
\end{lemma}
\begin{proof}
Induct on $i$, and for fixed $i$ on the finite cardinality $|S|$.  For $i=1$, formula \eqref{app:eq:initial-betti} gives $\beta_1^R(Q_S)=|S|<\infty$.  For $i\ge2$, the case $S=\varnothing$ has $Q_S=R$ and $\beta_i^R(Q_S)=0$.

Suppose $S\ne\varnothing$ and choose $v\in S$.  Set $J=S\setminus\{v\}$ and $T=J\cup N(v)$.  Then
\[
\mathcal B_{i-1}(J)\subseteq\mathcal B_{i-1}(S),
\qquad
\mathcal B_{i-2}(T)\subseteq\mathcal B_{i-1}(S).
\]
In particular $T$ is finite, since $i-1\ge1$.  The inner induction gives $\beta_i^R(Q_J)<\infty$, and the outer induction gives $\beta_{i-1}^R(Q_T)<\infty$.  Recurrence \eqref{app:eq:betti-recurrence} now yields $\beta_i^R(Q_S)<\infty$.
\end{proof}

Let
\[
R=\frac{\mathbb F[x_1,x_2,\ldots,y_1,y_2,\ldots]}
{(x_jx_{j+1},\ x_1y_1,\ y_1y_i\ :\ i,j\ge1)},
\qquad \mathfrak m=R_+,\qquad D=R_{\mathfrak m}.
\tag{8}\label{app:eq:example-rings}
\]
There are graph edges $x_j\! -\! x_{j+1}$, $x_1\! -\! y_1$, and $y_1\! -\! y_i$ for every $i\ge1$, including the loop at $y_1$.  No other edges occur.  In particular,
\[
N(x_1)=\{x_2,y_1\},\quad
N(x_r)=\{x_{r-1},x_{r+1}\}\ (r\ge2),\quad
N(y_1)=\{x_1,y_1,y_2,y_3,\ldots\}.
\tag{9}\label{app:eq:example-neighbors}
\]

\begin{proposition}\label{app:prop:betti}
For each $r\ge1$, put $S^{(0)}=N(x_r)$ and let
$M_r=Q_{S^{(0)}}=R/I_{N(x_r)}$.  Then
\[
\beta_j^R(M_r)<\infty\quad(0\le j\le r),
\qquad
\beta_{r+1}^R(M_r)=\aleph_0.
\tag{10}\label{app:eq:example-betti}
\]
\end{proposition}
\begin{proof}
The set $S^{(0)}$ is finite.  If $r=1$, it contains $y_1$, while each $y_i$ with $i\ge2$ lies at distance one.  If $r\ge2$, the shortest path from $S^{(0)}$ to $y_1$ is
\[
x_{r-1},x_{r-2},\ldots,x_1,y_1,
\]
with $r-1$ edges, and every $y_i$ with $i\ge2$ is at distance exactly $r$.  Since each $x_j$ has finitely many neighbors, the radius-$(r-1)$ ball $\mathcal B_{r-1}(S^{(0)})$ is finite.  Lemma~\ref{app:lem:radius} proves $\beta_r^R(M_r)<\infty$; its smaller-radius versions show $\beta_j^R(M_r)<\infty$ for $1\le j<r$, while $\beta_0^R(M_r)=1$.

For the first infinite Betti number, recursively form
\[
S^{(a)}=(S^{(a-1)}\setminus\{v_a\})\cup N(v_a).
\]
For $r\ge2$, choose the sequence of $r$ vertices
\[
v_1=x_{r-1},\ v_2=x_{r-2},\ \ldots,\ v_{r-1}=x_1,\ v_r=y_1;
\]
for $r=1$, choose only $v_1=y_1$.  At each step, $v_a\in S^{(a-1)}$: the preceding substitution introduces the next variable along the chain, and the substitution at $x_1$ introduces $y_1$.  All sets $S^{(a)}$ for $a<r$ are finite because each $N(x_j)$ is finite.  But $S^{(r)}$ contains the infinite set $N(y_1)$.

Applying \eqref{app:eq:betti-recurrence} successively and retaining its second summand gives
\[
\beta_{r+1}^R(M_r)\ge\beta_r^R(Q_{S^{(1)}})\ge\cdots
\ge\beta_1^R(Q_{S^{(r)}})=|S^{(r)}|=\aleph_0.
\]
The total set of monomials in $R$ is countable, hence every term of the bar complex in Lemma~\ref{app:lem:bar} is a countably generated $\mathbb F$-vector space.  Thus no $\Tor$ group in question has dimension greater than $\aleph_0$, and the final Betti number is exactly $\aleph_0$.
\end{proof}

\begin{lemma}\label{app:lem:minimal}
Let $H=\bigoplus_{j\ge0}H_j$ be a nonnegatively graded algebra over a field with $H_0=\mathbb F$, possibly with infinitely many generators.  Any graded $H$-module generated in degrees bounded below has a graded free resolution $F_\bullet$ whose matrices have entries in $H_+$.  For such a resolution,
\[
\Tor_i^H(\mathbb F,M)\cong\mathbb F\otimes_HF_i
\]
as graded vector spaces, so the number of homogeneous free summands of $F_i$ in each degree equals the corresponding graded Tor dimension.
\end{lemma}
\begin{proof}
The existence of minimal graded projective resolutions for bounded-below
graded modules over algebras with finite-dimensional degree-zero part
is established by Reyes--Rogalski
\cite[Section~2.1, preceding Lemma~2.2]{ReyesRogalski}.
Here $H_0=\mathbb F$, so the graded projective terms are graded free.
Minimality identifies their homogeneous summands with the graded Tor
groups; see \cite[Lemma~2.6 and its proof]{ReyesRogalski}.
\end{proof}

\begin{theorem}\label{app:thm:FP}
For $R$ and $D$ as in \eqref{app:eq:example-rings}, and every $r\ge1$,
\[
(x_r)D\in\FP_r(D)
\qquad\text{and}\qquad
(x_r)D\notin\FP_{r+1}(D).
\tag{11}\label{app:eq:finite-presentation}
\]
More precisely,
\[
\dim_\mathbb F\Tor_j^D(\mathbb F,(x_r)D)<\infty\quad(0\le j\le r),
\qquad
\dim_\mathbb F\Tor_{r+1}^D(\mathbb F,(x_r)D)=\aleph_0.
\tag{12}\label{app:eq:localized-betti}
\]
\end{theorem}
\begin{proof}
The same monomial-basis argument as in Lemma~\ref{app:lem:colon} gives
\[
\operatorname{ann}_R(x_r)=I_{N(x_r)}.
\]
Therefore multiplication by $x_r$ induces the graded isomorphism
\[
(x_r)\cong Q_{N(x_r)}(-1).
\tag{13}\label{app:eq:principal-ideal}
\]
Proposition~\ref{app:prop:betti} and Lemma~\ref{app:lem:minimal} imply that $(x_r)$ has a graded free resolution whose terms through homological degree $r$ have finite rank, while the term in degree $r+1$ necessarily has countably infinite rank.  In particular, $(x_r)\in\FP_r(R)$.

Localizing the first $r+1$ terms at $\mathfrak m$ gives a finite free partial resolution of $(x_r)D$, so $(x_r)D\in\FP_r(D)$.  By flat localization and the identification $R/\mathfrak m=\mathbb F$,
\[
\Tor_j^D(\mathbb F,(x_r)D)
\cong D\otimes_R\Tor_j^R(\mathbb F,(x_r))
\cong\Tor_j^R(\mathbb F,(x_r)).
\tag{14}\label{app:eq:localized-tor}
\]
The last isomorphism holds because each Tor group on the right is a $\mathbb F$-vector space, and every element of $R\setminus\mathfrak m$ acts invertibly on it.  Equation \eqref{app:eq:localized-betti} follows from Proposition~\ref{app:prop:betti}.

Finally, suppose $(x_r)D$ admitted a finite free partial resolution
through degree $r+1$ in the ungraded module category.  Extend that partial resolution to a (possibly infinitely generated) free resolution in higher degrees.  After tensoring with $\mathbb F$, the vector space $\Tor_{r+1}^D(\mathbb F,(x_r)D)$ would be a subquotient of $\mathbb F\otimes_D D^{a_{r+1}}$ for some finite $a_{r+1}$, hence finite-dimensional.  This contradicts \eqref{app:eq:localized-betti}.  Thus $(x_r)D$ fails $\FP_{r+1}$ in the ordinary, ungraded module
category, as asserted in \eqref{app:eq:finite-presentation}.
\end{proof}

\section*{Acknowledgements}
The research is supported by the Ministry of Science and Higher Education
of the Russian Federation (Goszadaniye) 075-00337-20-03, project
No.~075-00141-26-00.

\vspace{1.5em}
\noindent\textsc{Alexander Pavlov}\\
HSE University, Moscow\\
Moscow Institute of Physics and Technology\\
\href{mailto:abpavlov@hse.ru}{\nolinkurl{abpavlov@hse.ru}}

\medskip
\noindent\textsc{Mikhail Dusman}\\
Moscow Institute of Physics and Technology\\
\href{mailto:dusman.ma@phystech.edu}{\nolinkurl{dusman.ma@phystech.edu}}
\end{document}